\documentclass[12pt]{article}
\usepackage[francais]{babel}
\usepackage{amsmath,amsfonts,amssymb,amsthm}
\newcommand{\N}{\mathbb{N}}
\newcommand{\R}{\mathbb{R}}
\newcommand{\D}{\mathbb{D}}
\newcommand{\M}{\mathcal{M}}

\newcommand{\rg}{\rightarrow}

\begin{document}

\title{Un th\'eor\`eme de Helson\\ pour des s\'eries de Walsh}
\author{Jean--Pierre Kahane\\
Laboratoire de Math\'ematique,\\
Universit\'e Paris--Sud \`a Orsay}
\date{}
\maketitle

Henry Helson a \'etabli en 1954 le th\'eor\`eme suivant : \textit{Si les sommes partielles d'une s\'erie trigonom\'etrique sont positives, les coefficients tendent vers z\'ero }\cite{hel}.  C'est une des perles de la th\'eorie des s\'eries trigonom\'etriques, encore mise en valeur par le fait que l'hypoth\`ese n'entra\^{i}ne pas que la s\'erie trigonom\'etrique est une s\'erie de Fourier--Lebesgue (Katznelson 1965 \cite{kat}).

Nous allons montrer l'analogue des th\'eor\`emes de Helson et de Katznelson pour les s\'eries de Walsh, avec un compl\'ement au th\'eor\`eme de Katznelson qui s'\'etend aux s\'eries trigonom\'etriques.

\vskip2mm

\textsc{Th\'eor\`eme} I.--- \textit{Si les sommes partielles d'une s\'erie de Walsh sont positives, les coefficients tendent vers z\'ero.}

\vskip2mm

\textsc{Th\'eor\`eme} II.---  \textit{Soit $\psi$ une application croissante de $\R^+$  sur $\R^+$, telle que $\lim\limits_{x\rg 0} (\psi(x)/x^2)=0$. Alors a) il existe une mesure de probabilit\'e singuli\`ere $\mu$ sur $[0,1]$ dont les sommes partielles de la s\'erie de Walsh sont positives, et dont les coefficients de Fourier--Walsh, $\widehat\mu(n)$, v\'erifient $\sum\limits_0^\infty \psi(|\widehat\mu(n)|)<\infty$ b) le m\^eme \'enonc\'e vaut pour les s\'eries et les coefficients de Fourier au sens usuel.}

\vskip2mm

Pr\'ecisons les notations pour les s\'eries de Walsh. Au lieu de $[0,1]$, il est commode de les d\'efinir sur le groupe multiplicatif $\{-1,1\}^{\N^*}$, que nous d\'esignerons par $\D$. Les fonctions coordonn\'ees $r_1,r_2,\ldots$ sont les fonctions de Rademacher, et elles engendrent par multiplication les fonctions de Walsh $w_n$ $(n=0,1,2,\ldots)$ qui sont les caract\`eres de $\D$. On ordonne ainsi les $w_n=1$, $w_1=r_1$, $w_2=r_2$, $w_3=r_1r_2$, $w_4=r_3$ etc, c'est--\`a--dire
$$
w_n = \prod\ r_j^{\alpha_j}\,,\quad \alpha_j=0\ \hbox{ou}\ 1 \ \ \hbox{et}\ \ \sum \ \alpha_j <\infty\,. \leqno(1)
$$
Aux entiers  $n$ on associe ainsi les suites finies $(\alpha_j)$ de $0$  et de $1$, et l'ordre croissant des $n$ est l'ordre lexicographique inverse des mots $(\alpha_j)$.

Une s\'erie de Walsh est une s\'erie formelle de la forme $\sum\limits_0^\infty c_n w_n$, o\`u les coefficients $c_n$ sont r\'eels ou complexes. Nous nous bornerons aux $c_n$ r\'eels.  Les fonctions de Walsh op\'erant par multiplication sur les s\'eries de Walsh. Nous aurons besoin du lemme suivant.

\vskip2mm

\textsc{Lemme}.--- \textit{Soit $S$ une s\'erie de Walsh et $w_m$ une fonction de Walsh. Ecrivons $w_mS$ sous la forme
$$
w_m\, S= \sum_0^\infty d_n w_n\,.
$$
Alors, pour chaque $k$, la somme partielle d'ordre $2^k$, $\sum\limits_0^{2^k-1} d_n\, w_n$, est le produit par $w_m$ d'une diff\'erence de sommes partielles de la s\'erie~$S$.}

\vskip2mm

\textbf{\textit{Preuve}}. La s\'erie $\sum\limits_0^\infty d_n\,w_n\,w_m$ n'est autre que la s\'erie $S$ dont on a modifi\'e l'ordre des termes. Reste \`a montrer que les $w_n\,w_m$ $(n=0,1,\ldots,2^k-1)$ constituent, \`a l'ordre pr\`es, une suite de $2^k$ fonctions de Walsh cons\'ecutives. Pour cela, \'ecrivons les $w_n$ sous la forme (1) (avec ici $\alpha_j=0$ ou $1$ pour $j\le k$ et $\alpha_j=0$ pour $j>k$) et $w_m = \prod r_j^{\beta_j}$ $(\beta_i =0$ ou 1).  Les $w_n\,w_m$ s'\'ecrivent $\prod r_j^{\gamma_j}$ avec $\gamma_j =0$ ou $1$ et $\gamma_j = \alpha_j+\beta_j$ modulo $1$. Ainsi $(\gamma_j)_{j=1,2,\ldots,k}$ parcourt $\{0,1\}^k$ quand $(\alpha_j)_{j=1,2,,\ldots,k}$ parcourt $\{0,1\}^k$, tandis que $(\gamma_j)_{j>k}= (\beta_j)_{j>k}$. Donc les indices des $w_n\,w_m$, qui s'\'ecrivent $\sum\limits_1^k +\sum\limits_{k+1}^\infty\,\gamma_j\, 2^{j-1}$, parcourent un segment des entiers de longueur $2^k$, CQFD.

A toute s\'erie de Walsh $S$ est associ\'ee une martingale dyadique constitu\'ee par ses sommes partielles d'ordre $2^k$ $(k=0,1,2,\ldots)$
$$
M_k =M_k (S) = \sum_0^{2^k-1} c_n\,w_n\,,
$$
et on obtient ainsi toutes les martingales dyadiques d\'efinies sur $\D$. Rappelons des propri\'et\'es des martingales dyadiques dont nous nous servirons.

\vskip2mm

\textbf{P1.}\  \ $S$ est la s\'erie de Fourier--Walsh d'une mesure de Radon r\'eelle  sur $\D$, c'est--\`a--dire $c_n=\int w_n d\mu$, $\mu\in \M(\D)$, si et seulement si les $M_k$ sont born\'es dans $L^1(\D)$. Dans ce cas, les $M_k$ tendent presque partout sur $\D$ vers la densit\'e de la partie absolument continue de $\mu$, et les $M_k$ tendent vers $\mu$ dans $\M(\D)$ au sens faible, comme formes lin\'eaires sur $C(\D)$.

\vskip2mm

\textbf{P2.}\  \ $S$ est la s\'erie de Fourier--Walsh d'une fonction r\'eelle int\'egrable sur $\D$, soit $c_n=\int f\, w_n$, $f\in L^1(\D)$, si et seulement si les $M_k$ sont uniform\'ement int\'egrables. Dans ce cas, les $M_k$ tendent vers $f$ presque partout et dans $L^1(\D)$.

\vskip2mm

\textbf{P3.}\ \ $S$ est la s\'erie  de Fourier--Walsh d'une mesure positive $\mu\in \M^+(\D)$ si et seulement si les $M_k$ sont positives, $M_k\ge 0$.

\vskip2mm

Ici comme dans la suite, positif signifie $\ge0$.

\vskip2mm

\textbf{\textit{Preuve du th\'eor\`eme }I}.

Supposons les sommes partielles de la s\'erie $S$ positives, et de plus $c_0=1$. Alors $S$ est la s\'erie de Fourier--Walsh d'une mesure de probabilit\'e $\mu\in \M_1^+(\D)$, soit $c_n=\widehat \mu(n)$, et
$$
M_k = M_k(r_1,r_2,\ldots,r_k) = \sum_0^{2^k-1} \widehat\mu(n) w_n\,.
$$

Ecrivons
\[
\begin{array}{rl}
  M_{k+1} &=  M_k +r_{k+1}N_k \,, \ \ N_k = N_k(r_1,r_2,\ldots r_k) \\
  N_k&= \displaystyle\sum_0^{2^k-1} \widehat N_k(m) w_m   \\
  N_k^*&= \sup\limits_{0\leq n<2k} \Big| \displaystyle \sum_0^n \widehat N_k(m) w_m\Big|   
\end{array}
\]
L'hypoth\`ese que les sommes partielles de $S$ soient positives se traduit par
$$
N_k^* \le M_k \qquad (k=0,1,\ldots) \leqno(2)
$$

Supposons que les $\widehat \mu(n)$ ne tendent pas vers $0$, c'est--\`a--dire qu'il existe un $a>0$, une suite  strictement croissante d'entiers $k_j$, et des entiers $n_j\in [2^{k_j}, 2^{k_j+1}[$ tels que $|\widehat \mu(n_j)|\ge a$, et tentons d'\'etablir une contradiction.

Supposons d'abord $n_j=2^{k_j}$. Les $N_{k_i}$ sont born\'es dans $L^1(\D)$ et $|\widehat N_{k_j}(0)|\ge a$ puisque $\widehat N_{k_i}(0)=\widehat\mu(2^{k_j})$. Quitte \`a remplacer la suite $(k_j)$ par une sous--suite, nous pouvons supposer que les $N_{k_j}$ convergent faiblement vers une mesure $\nu\in \M(\D)$. Ainsi
$$
\widehat \nu(n)=\lim_{i\rg\infty} \widehat N_{k_i}(n)
$$
et en particulier $\widehat \nu(0)\neq 0$.

Les sommes
$$
\nu_k= \sum_0^{2^k-1} \widehat\nu(n)w_n\qquad (k=0,1,2,\ldots)
$$
forment une martingale dyadique born\'ee dans $L^1(\D)$ (propri\'et\'e P1) et,  pour chaque $k$,
$$
\nu_k=\lim_{j\rg\infty} \sum_0^{2^{k-1}} \widehat N_{k_i}(n) w_n\,.
$$
L'hypoth\`ese de positivit\'e, sous la forme (2), entra\^{i}ne
$$
\Big | \sum_0^{2^k-1} \widehat N_{k_i}(n) w_n\Big | \le M_{k_i}
\leqno(3)
$$
lorsque $k_j\ge k$, donc $|\nu_k|\le M_{k_i}$.

Or les $M_{k_i}$ convergent presque partout vers une $f\in L^1$, la densit\'e de la partie absolument continue de $\mu$ (propri\'et\'e P1). Donc
$$
|\nu_k| \le f\,.
$$
Cela entra\^{i}ne que les $\nu_k$ sont uniform\'ement int\'egrables, donc convergent dans $L^1(\D)$ (propri\'et\'e P2), donc que $\nu$ est absolument continue.

D'autre part, si l'on d\'ecompose $\mu$ en sa partie absolument continue et sa partie singuli\`ere, $\mu=\mu^a+\mu^s$, on peut \'ecrire \'egalement $M_k=M_k^a +M_k^s$ et $N_k= N_k^a+ N_k^s$. Les $M_k^a$ convergent vers $f$ dans $L^1(\D)$, donc les $N_k^a$ tendent vers $0$ dans $L^1(\D)$, donc $\nu$ est la limite faible des $N_{k_j}^s$. Or, pour tout $\varepsilon>0$ et tout $k$,
$$
(N_k^s >\varepsilon) \subset \big(M_k^s > {\varepsilon\over2}\big) \cup \big(M_{k+1}^1 > {\varepsilon\over2}\big)
$$
et la mesure du second membre tend vers 0 quand $k\rg \infty$. Donc $\nu$ est singuli\`ere. La contradiction est \'etablie dans le cas particulier $n_j = 2^{k_i}$.

Passons au cas g\'en\'eral. On consid\`ere maintenant les 
$$
N_{k_j}' = w_{m_j} N_{k_j}\,,\qquad m_j = n_j- 2^{k_j}\,.
$$
Ainsi $\widehat N_{k_j}'(0) = \widehat N_{k_j} (m_j) = \widehat \mu(n_j)$. Comme ci--dessus, quitte \`a restreindre la suite $(n_j)$, les $N_{k_j}'$ convergent faiblement vers une mesure $\nu'$ singuli\`ere, non nulle puisque $\widehat \nu'(0)\neq 0$. Pour montrer que $\nu'$ est absolument continue, le point crucial est l'analogue de (3) que l'on obtient en appliquant le lemme \`a $N_{k_j}$ (pour $S$) et $w_{m_j}$ (pour $w_m$). On obtient ainsi, pour $k\ge k_j$,
$$
\bigg|\sum_0^{2^k-1} \widehat N_{k_j}'(n) w_n \bigg| \le 2\ M_{k_j}\,,
\leqno(4)
$$
ce qui permet d'achever la d\'emonstration que $\nu'$ est absolument continue. La contradiction est ainsi \'etablie dans le cas g\'en\'eral, et cela ach\`eve la preuve du th\'eor\`eme~1.

\vskip2mm

\textbf{Remarque.}\ \ Cette preuve est calqu\'ee sur celle de Helson. Comme ici, Helson met en \'evidence, sous l'hypoth\`ese que les coefficients ne tendent pas vers z\'ero, une mesure $\nu$ non nulle qui est \`a la fois singuli\`ere et absolument continue. Mais la d\'emonstration de la continuit\'e absolue est diff\'erente. Helson utilise un th\'eor\`eme de Fr\'ed\'eric et Marcel Riesz qui appartient \`a la th\'eorie des fonctions analytiques. On utilise ici la th\'eorie des martingales. L'emploi en parall\`ele des martingales et des fonctions analytiques est classique en analyse harmonique depuis les th\'eor\`emes de Paley et de Littlewood--Paley.

Une autre diff\'erence entre les preuves est l'utilisation de l'hypoth\`ese. La positivit\'e  des sommes partielles entra\^{i}ne que les sommes partielles sont born\'ees dans $L^1$, et Helson n'utilise rien d'autre. Au contraire,  nous avons utilis\'e de mani\`ere essentielle une autre cons\'equence de la positivit\'e, la formule (2) (\'equivalente \`a la positivit\'e des sommes partielles). 

\vskip2mm

\textbf{\textit{Preuve du th\'eor\`eme }II}.

\vskip1mm

\textbf{\textit{Partie a)}}

On construit la mesure $\mu$ comme produit infini
$$
\mu = \prod (1+X_k)
$$
o\`u les $X_k$ sont de polyn\^omes de Walsh ind\'ependants \`a valeurs dans l'intervalle $[-1,1]$ et de valeur moyenne nulle : $|X_k|\le 1$ et $EX_k =0$. Posons $EX_k^2 =\sigma_k^2$. Lorsque $\sum \sigma_k^2 <\infty$, le produit infini converge dans $L^2(\D)$, donc $\mu$ est absolument continue. Lorsque $\sum \sigma_k^2=\infty$, $\mu$ est une mesure de probabilit\'e singuli\`ere.

Voici une d\'emonstration rapide de ce dernier fait, tir\'ee de \cite{pey}. Supposons $\sum \sigma_k^2=\infty$. Les ${X_k\over\sigma_k}$ forment un syst\`eme orthonormal dans $L^2(\D,\lambda)$, o\`u $\lambda$ d\'esigne la mesure de Haar sur $\D$. On v\'erifie que les ${X_k\over\sigma_k} -\sigma_k$ sont de valeur moyenne nulle et deux \`a deux orthogonales dans $L^2(\D,\mu)$ :
$$
E_\mu\Big({X_k\over\sigma_k} - \sigma_k\Big) = E_\lambda\Big(\Big({X_k\over\sigma_\kappa} - \sigma_k\Big) (1+X_k)\Big)=0
$$
et, pour $k\neq k'$,
\[\!\!
\begin{array}{rl}
 E_\mu\big(\big({X_k\over\sigma_k}\!-\!\sigma_k\!\big)\big({X_{k'}\over\sigma_{k'}}\! -\! \sigma_{k'}\big)\big) \kern-3mm&=\!E_\lambda\big(\big({X_k\over\sigma_k} -\sigma_k \big)\big({X_{k'}\over \sigma_{k'}} - \sigma_{k'}\big)(1+X_k)(1+X_{k'})\big)  \\
  &=\!E_\lambda\big(\big({X_k\over\sigma_k}\! -\!\sigma_k \big)(1\!+\!X_k)\big) E_\lambda \big(\big({X_{k'}\over \sigma_{k'}}\! -\! \sigma_{k'}\big)(1\!+\!X_{k'})\big)\!=\!0. \\  
\end{array}
\]
De plus,
\[
\begin{array}{c}
E_\mu\big(\big({X_k\over\sigma_k} - \sigma_k\big)^2\big) = E_\lambda\big(\big( {X_k\over\sigma_k}-\sigma_k\big)^2 (1+X_k)\big)   \\
\le 2\ E_\lambda\big({X_k\over\sigma_k} - \sigma_k\big)^2 = 2(1+ \sigma_k^2) \le 4\,.    \\  
\end{array}
\]
Pour toute suite $(b_k) \in \ell^2$, les s\'erie $\sum b_k{X_k\over\sigma_k}$ et $\sum b_k \big({X_k\over\sigma_k} - \sigma_k\big)$ convergent respectivement dans $L^2(\D,\lambda)$ et dans $L^2(\D,\mu)$. Si $\mu$ n'\'etait pas orthogonale \`a $\lambda$, il existerait un point de $\D$ et une suite d'entiers $n_i$ tels que les sommes partielles d'ordre $n_i$ des deux s\'eries convergent en ce point. Par diff\'erence, les sommes partielles d'ordre $n_i$ de la s\'erie $\sum b_k\, \sigma_k$ convergeraient. Or on peut choisir les $b_k >0$, $(b_k)\in \ell^2$, de fa\c con que la s\'erie $\sum b_k\,\sigma_k$ diverge. Donc $\mu \bot \lambda$.

Remarquons que l'hypoth\`ese d'ind\'ependance des $X_k$ peut \^etre remplac\'ee par une condition d'orthogonalit\'e forte, \`a savoir
$$
\int \prod X_k^{\alpha_k} = 0
$$
pour toutes les suites $(\alpha_k)$ constitu\'ees de $0$, $1$ et 2, et finies $(\sum \alpha_k < \infty)$, contenant au moins un 1 et au plus deux 2. C'est sous cette  forme que la m\'ethode a \'et\'e introduite et utilis\'ee par Peyri\`ere dans l'\'etude de la singularit\'e mutuelle des produits de Riesz \cite{pey}.

On imposera donc la condition $\sum \sigma_k^2 =\infty$, et il s'agit maintenant de construire les $X_k$ de fa\c con que les sommes partielles de la s\'erie de Fourier--Walsh de $\mu$ soient positives, et que les coefficients v\'erifient la majoration $|\widehat\mu(n)| \le \psi(n)$. Pour cela, il est commode de disposer de polyn\^omes de Walsh de la forme suivante :
$$
\varphi=\varphi(r_1,r_2,\ldots r_\ell) = \sum_1^{2^\ell}\varepsilon_n w_n\,, \quad \varepsilon_n=\pm1
$$
pour lesquels
$$
\| \varphi\|_v = \sup_{1\le p\le 2^\ell} \|\sum_1^p \varepsilon_n w_n\|_\infty < C\, 2 ^{\ell/2} = C\|\varphi\|_2\,,
$$
$C$ \'etant une constante absolue. Voici une construction classique de tels polyn\^omes : on pose $P_0 =Q_0=1$, $P_{\ell+1} = P_\ell +r_{\ell+1} Q_\ell$, $Q_{\ell+1} = P_\ell - r_{\ell+1}Q_\ell$ $(2=0,1,\ldots)$, on v\'erifie que $P_\ell^2 +Q_\ell^2 = 2^{\ell+1}$, d'o\`u
$$
\|P_\ell\|_v \le 2^{\ell/2}+ 2^{(\ell-1)/2} +\cdots \le 2^{\ell/2} {\sqrt{2}\over\sqrt{2}-1}
$$
et on pose $\varphi(r_1,r_2,\ldots, r_\ell)=r_1P_\ell$. C'est la construction donn\'ee dans le cas trigonom\'etrique par Harold Shapiro puis par Walter Rudin, et la suite $(\varepsilon_r)$ s'appelle la suite de Rudin--Shapiro.

Etant donn\'e un ensemble d'entiers positifs $J$, de cardinal $2^\ell$, on d\'esignera par $\varphi((r_j),\, j\in J)$ le polyn\^ome de Walsh obtenu \`a partir de $\varphi(r_1,\, r_2,\ldots r_\ell)$ en substituant \`a $r_1,\,r_2,\ldots r_\ell$ les $r_j$, $j\in J$, dans l'ordre croissant des $j$. La construction de $\mu$ va d\'ependre essentiellement du choix d'une suite tr\`es rapidement croissante d'entiers $\ell_k$, que nous ferons plus tard. Pour chaque $k$, soit $J_k$ un ensemble d'entiers, de cardinal $2^{\ell_k}$, situ\'e \`a droite de $J_{k-1} : \inf J_k>\sup J_{k-1}$. Posons
\[
\begin{array}{l}
\displaystyle  a_k = {1\over2C} 2^{-\ell_k/2}   \\
X_k = a_k \ \varphi((r_j),\ j\in J_k)\,.   \\  
\end{array}
\]
Explicitons les normes de $X_k$ dans $L^2$, dans $U$ (maximum des valeurs absolues des sommes partielles), dans $A$ (somme des valeurs absolues des coefficients) et dans $PM$ (sup des valeurs absolues des coefficients) :
\[
\begin{array}{ll}
  \|X_k\|_2     & = \displaystyle {1\over2C}\, ,    \\
  \noalign{\vskip2mm}
   \|X_k\|_U     &< \displaystyle {1\over2}\, ,   \\
     \noalign{\vskip2mm}
   \|X_k\|_A     & = \displaystyle {1\over2C} 2^{\ell_k/2}\, , \\ 
     \noalign{\vskip2mm}
    \|X_k\|_{PM}     & = \displaystyle {1\over2C} 2^{-\ell_k/2}\, .
\end{array}
\]
Comme les $J_k$ sont disjoints, les $X_k$ sont bien ind\'ependants, et on a bien $EX_k=0$ et $\sum \sigma_k^2=\infty$.

Posons 
$$
\Pi_k= (1+X_1)(1+X_2)\cdots (1+X_k)\,.
$$
Une somme partielle de la s\'erie de Fourier--Walsh de $\mu$ dont l'ordre est compris entre $\sup J_j$ et $\sup J_{j+1}$ est la somme de $\Pi_k$ et d'une somme partielle de $\Pi_k X_{k+1}$, ce que nous \'ecrivons
$$
S_\cdot (\mu) = S_\cdot(X_{k+1}) = \Pi_k +S_{\cdot\cdot} (\Pi_k X_{k+1})\,.
$$
Cette derni\`ere somme partielle se d\'ecompose \`a son tour en
$$
S_{\cdot\cdot} (\Pi_k X_{k+1}) = \Pi_k S_{\cdot\cdot\cdot}(X_{k+1}) + S_{\cdot\cdot\cdot\cdot}(\Pi_k) \times \hbox{un coefficient de }X_{k+1}\,.
$$
Donc
\[
\begin{array}{rl}
 S_{\cdot}(X_{k+1})  & \ge \Pi_k -\Pi_k \|X_{k+1}\|_U - S_{\cdot\cdot\cdot\cdot}(\Pi_k)\|X_{k+1}\|_{PM}   \\
 \noalign{\vskip2mm}
  &\ge \displaystyle{1\over2} \Pi_k - {1\over2C} 2^{-(\ell_{k+1}/2)} \|\Pi_k\|_A\,.   \\
\end{array}
\]
Imposons la condition que, pour tout $k$,
$$
{1\over2C} \ 2^{-(\ell_{k+1}/2)} \|\Pi_k\|_A \le {1\over4} \inf \Pi_k\,. \leqno(5)
$$
Il en r\'esulte que $S_\cdot (X_{k+1}) \ge {1\over4} \Pi_k$, donc $S_\cdot(\mu)\ge 0$.

Ecrivons $\psi(x) =x^2\,\varepsilon(x)$. Alors
$$
\sum \varphi(\widehat\mu(n)) \le \sum_k \|\Pi_k \, X_{k+1}\|_2^2 \ \varepsilon(\|\Pi_k\, X_{k+1}\|_{PM})\,.
$$
Or
$$
\|\Pi_k\, X_{k+1}\|_2^2 \le {1\over4C^2} \|\Pi_k\|_A^2
$$
et
$$
\varepsilon(\|\Pi_k \, X_{k+1}\|)_{PM} \le \varepsilon (\|X_{k+1}\|_{PM}) = \varepsilon\Big({1\over2C} 2^{-(\ell_{k+1}/2)}\Big)\,.
$$
Si,  outre (5), on impose \`a la suite $(\ell_k)$ la condition
$$
\sum \|\Pi_k\|_A^2\ \varepsilon\Big({1\over2C}\ 2^{-(\ell_{k+1}/2)}\Big) <\infty\,,
\leqno(6)
$$
ce qui est possible, la conclusion du th\'eor\`eme est v\'erifi\'ee.

La conclusion de la partie a) du th\'eor\`eme est v\'erifi\'ee.
\vskip2mm

\textbf{\textit{Partie b)}}

\vskip1mm

On choisit ici des polyn\^omes trigonom\'etriques
$$
\varphi_\ell (t) = \sum_a^\ell \varepsilon_n \cos nt\,, \qquad \varepsilon_n=\pm1\,,
$$
avec la propri\'et\'e de
$$
\|\varphi_\ell\|_v \le C\ \ell^{1/2} = C \|\varphi_\ell\|_2\,,
$$
et on choisit
$$
X_k(t) = a_k \ \varphi_{\ell_k} (\ell_k t)\,,\qquad a_k = {1\over4C}\ 2^{-\ell_k/2}\,.
$$
Les $X_k$ ne sont plus des fonctions ind\'ependantes, mais, si la suite $(\ell_k)$ est assez rapidement croissante, elles sont orthogonales au sens fort qui a \'et\'e d\'ecrit ci--dessus. On d\'efinit donc
$$
\mu = \prod_1^\infty (1+a_k\ \varphi_{\ell_k}(\ell_k\, t))
$$
et on v\'erifie par les m\^emes calculs que ci--dessus que $\mu$ est singuli\`ere, que ses sommes partielles sont positives, et que $\sum \varphi(|\widehat\mu(n)|) <\infty$.

C'est exactement la m\'ethode de Katznelson.

Cet article a \'et\'e \'ecrit en septembre 2004 pour f\^eter les 50 ans du th\'eor\`eme de Helson et les 70~ans d'Yitzhak Katznelson (14 novembre 2004).

 \eject
 
 \vglue 2cm

\vskip4mm

\hfill\begin{minipage}{6,5cm}
Jean--Pierre Kahane

Laboratoire de Math\'ematique

Universit\'e Paris--Sud, B\^at. 425

91405 Orsay Cedex

\textsf{Jean-Pierre.Kahane@math.u-psud.fr}

\end{minipage}


\begin{thebibliography}{2}

\bibitem[1]{hel}
\textsc{Helson}, Henry  \textit{Proof of a conjecture of Steinhaus}, Proc. Nat. Acad. USA 40 (1954), 205--206.

\bibitem[2]{kat}
\textsc{Katznelson}, Yitzhak  \textit{Trigonometric series with positive partial sums}, Bull. Amer. Math. Soc. 71 (1965) 718--719.

\bibitem[3]{pey}
 \textsc{Peyri\`ere}, Jacques \textit{Etude de quelques propri\'et\'es des produits de Riesz,} Annales de l'Institut Fourier 25, 2 (1975), 127--169.
 
 
 



\end{thebibliography}
\end{document}